\documentclass[a4paper]{amsart}
\usepackage{amsmath}
\usepackage{booktabs}
\usepackage{amssymb}

\usepackage[small,nohug,heads=LaTeX]{diagrams}
\usepackage{hyperref}
\diagramstyle[labelstyle=\scriptstyle]

\newtheorem{MainThm}{Theorem}

\newtheorem{thm}{Theorem}[section]
\newtheorem{cor}[thm]{Corollary}
\newtheorem{lem}[thm]{Lemma}
\newtheorem{prop}[thm]{Proposition}
\theoremstyle{definition}
\newtheorem{defn}[thm]{Definition}
\theoremstyle{remark}
\newtheorem{rem}[thm]{Remark}

\newtheorem{example}[thm]{Example}
\numberwithin{equation}{section}

\newcommand{\bF}{\mathbb{F}}

\newcommand{\bP}{\mathbb{P}}

\newcommand{\bR}{\mathbb{R}}
\newcommand{\bS}{\mathbb{S}}

\newcommand{\bZ}{\mathbb{Z}}

\newcommand{\MT}[2]{\bold{MT #1}(#2)}

\newcommand\lra{\longrightarrow}
\newcommand\trf{\mathrm{trf}}
\newcommand\Diff{\mathrm{Diff}}
\newcommand\Th{\mathrm{Th}}
\newcommand\colim{\mathrm{colim \,}}
\newcommand\hocolim{\mathrm{hocolim \,}}

\newcommand\Ker{\mathrm{Ker}}
\newcommand\Cotor{\mathrm{Cotor}}
\newcommand\Rk{\mathrm{Rank}}
\newcommand{\cotensor}{\Box}

\newcommand{\hcoker}{/\!\!/}
\newcommand{\hker}{\backslash\!\!\backslash}

\title{The homology of the stable non-orientable mapping class group}
\author{Oscar Randal-Williams}
\thanks{Supported by an EPSRC Studentship, DTA grant number EP/P502667/1}
\email{randal-w@maths.ox.ac.uk}

\address{Mathematical Institute\\
	24-29 St Giles'\\
	Oxford\\
	OX1 3LB\\
	United Kingdom}

\begin{document}
\begin{abstract}
Combining results of Wahl, Galatius--Madsen--Tillmann--Weiss and Korkmaz one can identify the homotopy-type of the classifying space of the stable non-orientable mapping class group $\mathcal{N}_\infty$ (after plus-construction). At odd primes $p$, the $\bF_p$-homology coincides with that of $Q_0(\mathbb{H \, P}^\infty_+)$, but at the prime 2 the result is less clear. We identify the $\bF_2$-homology as a Hopf algebra in terms of the homology of well-known spaces. As an application we tabulate the integral stable homology of $\mathcal{N}_\infty$ in degrees up to six.

As in the oriented case, not all of these cohomology classes have a geometric interpretation. We determine a polynomial subalgebra of $H^*(\mathcal{N}_\infty ; \bF_2)$ consisting of geometrically-defined characteristic classes.
\end{abstract}
\maketitle

\section{Introduction}

The mapping class groups $\mathcal{N}_g$ of non-orientable surfaces are not as widely studied as their counterparts for oriented surfaces, but with Wahl's proof \cite{Wahl} of homological stability for these groups one can apply the machinery of Madsen and Weiss \cite{MW} used to prove the Mumford conjecture, or its more concise variant \cite{GMTW}, to study their stable homology. Together these results show that the homology of $\mathcal{N}_\infty$ coincides with that of a component of an infinite loop space, $\Omega_0^\infty \MT{O}{2}$, which we define in Section \ref{MTODef}.

Inspired by Galatius' calculation \cite{galatius-2004} of the homology of the stable oriented mapping class group $\Gamma_\infty$, we calculate the mod 2 homology of the stable non-orientable mapping class group $\mathcal{N}_\infty$. The odd-primary homology is much simpler: it coincides with that of $Q_0(\mathbb{H} \mathbb{P}^\infty_+)$, and is discussed in Section \ref{IntegralDiscussion}. Here $Q(X) := \colim \Omega^n \Sigma^n X$ denotes the free infinite loop space on the space $X$, the subscript 0 denotes the connected component of the constant loop, and the subscript $+$ denotes the addition of a disjoint basepoint. The homology $H_*(Q(X) ; \bF_p)$ is completely known \cite{CLM} as a functor of $H_*(X ; \bF_p)$.

We adopt the convention that all (co)homology groups are with $\bF_2$ coefficients unless otherwise stated. Our main theorems give a determination of the $\bF_2$-homology of $\Omega_0^\infty \MT{O}{2}$ as a Hopf algebra:
\begin{MainThm}\label{ExactSeqMTO1}
There is an exact sequence of Hopf algebras
$$\bF_2 \lra H_*(\Omega_0^\infty \MT{O}{1}) \overset{\bar{\omega}_*} \lra H_*(Q_0(BO_{1 +})) \overset{\bar{\partial}_*} \lra H_*(Q_0(S^0)) \lra \bF_2$$
which is split as algebras.
\end{MainThm}

\begin{MainThm}\label{ExactSeqMTO2}
There is an exact sequence of Hopf algebras
$$\bF_2 \lra H_*(\Omega_0^\infty \MT{O}{2}) \overset{\omega_*} \lra H_*(Q_0(BO_{2 +})) \overset{\partial_*} \lra H_*(\Omega_0^\infty \MT{O}{1}) \lra \bF_2$$
which is split as algebras.
\end{MainThm}

We also give formulae for the maps $\bar{\partial}_*$ and $\partial_*$. Although these theorems determine the required homology groups, they do not do so in a very explicit manner. Thus we tabulate the first six (co)homology groups:

\vspace{0.5cm}

\begin{center}
\begin{tabular}{lcccccc}
\hline
Degree $i$ & 1 & 2 & 3 & 4 & 5 & 6 \\ \toprule
$\Rk QH_i(\Omega_0^\infty \MT{O}{2} ; \bF_2)$ & 1 & 2 & 3 & 3 & 5 & 6 \\ 
$\Rk H_i(\Omega_0^\infty \MT{O}{2} ; \bF_2)$ & 1 & 3 & 6 & 12 & 23 & 45 \\
\midrule
$H_i(\Omega_0^\infty \MT{O}{2} ; \bZ)$ & $\bZ_2$ & $\bZ_2^2$ & $\bZ_2^3 \oplus \bZ_{12}$ & $\bZ_2^7 \oplus \bZ$ & $\bZ_2^{16}$ & $\bZ_2^{29}$ \\
$H^i(\Omega_0^\infty \MT{O}{2} ; \bZ)$ & 0 & $\bZ_2$ & $\bZ_2^2$ & $\bZ_2^3 \oplus \bZ_{12} \oplus \bZ$ & $\bZ_2^7$ & $\bZ_2^{16}$ \\ \bottomrule
\end{tabular}
\end{center}

\vspace{0.5cm}

The integral homology can be calculated via the Bockstein spectral sequence. Only the primes 2 and 3 appear in this range: for the prime 3 the homology coincides with that of $Q(\mathbb{H \, P}^\infty_+)$ and the Bockstein spectral sequence can be calculated from the known one for $\mathbb{H \, P}^\infty$, which is degenerate from the first page. For the prime 2 we compare the Bockstein spectral sequence for $Q(BO_{2 +})$ with the representing elements for classes belonging to $\Omega_0^\infty \MT{O}{2}$. We find that it degenerates at the $E^3$ page, giving the claimed homology groups.

We also study a family of mod 2 characteristic classes of surface bundles, $\chi_i$, defined as follows. For a surface bundle $F \lra E \overset{\pi} \lra B$ there is an associated vector bundle $V$ of first real cohomologies, with fibre $H^1(F_b ; \bR)$ over $b \in B$. We can define
$$\chi_i(E) := w_i (V) \in H^i(B ; \bF_2)$$
the $i^{th}$ Stiefel--Whitney class of $V$. The main theorem is that these are stably independent for the non-orientable mapping class groups, so account for one indecomposable generator of the cohomology of $B \mathcal{N}_\infty$ in each dimension.
\begin{MainThm}\label{Generators}
In the non-orientable mapping class groups, the map
$$\bF_2 [ \chi_1, \chi_2, \chi_3, ... ] \lra H^*(B\mathcal{N}_g ; \bF_2)$$
is an injection in the stable range $* \leq (g-3)/4$.
\end{MainThm}

\subsection{Acknowledgements}
I gratefully acknowledge the support of an EPSRC Studentship, DTA grant number
EP/P502667/1. I would like to thank the anonymous referee for their many helpful comments and suggestions.

\section{Recollections}

\subsection{Homological stability}\label{HomStab}
Write $N_{g, b}$ for the non-orientable surface consisting of the connected sum of $g$ copies of $\mathbb{RP}^2$ with $b$ discs removed. We call $g$ the genus of the surface, and $b$ the number of boundary components, which we suppress in notation if it is 0. The mapping class group of $N_{g, b}$ is the group of components of the topological group of self-diffeomorphisms that fix the boundary pointwise
$$\mathcal{N}_{g, b} : = \pi_0 (\Diff(N_{g,b}, \partial N_{g,b})).$$
As long as $b \geq 1$ there are two stabilisation maps between these groups, obtained by gluing either a $N_{1, 2}$ or a $N_{0, 3}$ to a boundary component and extending diffeomorphisms by the identity to the new surface
\begin{align*}
    \alpha : \mathcal{N}_{g, b} \lra \mathcal{N}_{g+1, b}\\
	 \beta :  \mathcal{N}_{g, b} \lra \mathcal{N}_{g, b+1}
\end{align*}
and $\beta$ has a left inverse obtained by gluing a disc over a boundary component
$$\delta : \mathcal{N}_{g, b+1} \lra \mathcal{N}_{g, b}.$$
We quote below the theorem of Wahl on the effect on these maps on group homology.
\begin{thm}[Theorem A of \cite{Wahl}]
For $b \geq 1$, on homology
\begin{enumerate}
    \item $\alpha$ gives a surjection in degrees $* \leq g/4$ and an isomorphism in degrees $* \leq (g-3)/4$.
	\item $\beta$ gives an isomorphism in degrees $* \leq (g-3)/4$.
	\item $\delta$ gives a surjection in degrees $* \leq (g-1)/4$ and an isomorphism in degrees $* \leq (g-5)/4$.
\end{enumerate}
\end{thm}
We will say ``the stable range" to indicate the range $* \leq (g-3)/4$ in which $\alpha$ induces an isomorphism.

\subsection{The homotopy-type of the stable non-orientable mapping class group}\label{MTODef}
Once we have homological stability for $\mathcal{N}_g$, the machinery of Galatius, Madsen, Tillmann and Weiss \cite{GMTW} identifies the stable homology of these groups.

In order to describe the result we must first introduce the spectrum $\MT{O}{2}$. Let $U_n^\perp \lra Gr_2(\mathbb{R}^{n+2})$ be the $n$-dimensional complement to the tautological 2-plane bundle $U_2$. The inclusion $Gr_2(\mathbb{R}^{n+2}) \lra Gr_2(\mathbb{R}^{n+1+2})$ pulls back $U_{n+1}^\perp$ to $\epsilon^1 \oplus U_n^\perp$, so there is an induced map on Thom spaces
$$Th(\epsilon^1 \oplus U_n^\perp) \cong S^1 \wedge Th(U_n^\perp) \lra Th(U^\perp_{n+1}).$$
This defines the spectrum $\MT{O}{2}$ with $(n+2)$nd space $Th(U_n^\perp)$. It is not hard to check that $\pi_0(\MT{O}{2}) \cong \bZ$, using the cofibre sequences of spectra (\ref{MTO1Seq} and \ref{MTO2Seq}) in the next section.

Let $N_g \lra E \lra B$ be a smooth fibre bundle with fibre $N_g$ a non-orientable surface of genus $g$. There is an embedding of $E$ in $\mathbb{R}^{n+2} \times B$ over $B$, with a tubular neighbourhood homeomorphic to $N_vE$, the vertical normal bundle of $E$. A Pontrjagin--Thom collapse gives a map
$$S^{n+2} \wedge B_+ \lra Th(N_v E)$$
and the vertical normal bundle can be classified by a map $f: E \lra G_2(\mathbb{R}^{n+2})$ such that $f^*U_n^\perp \cong N_vE$. Thus $f^* U_2$ is $T_v E$, the vertical tangent bundle. Composing with the classifying map we have a map
$$S^{n+2} \wedge B_+ \lra Th(U_n^\perp)$$
with adjoint $B \lra \Omega^{n+2} Th(U_n^\perp)$. Following this with the stabilisation map to $\Omega^\infty \MT{O}{2}$ defines the Madsen--Tillmann map
$$\alpha_E: B \lra \Omega^\infty \MT{O}{2}.$$
Applying this to the universal $F_g$ bundle over $B \mathcal{N}_g$ gives
$$\alpha_g : B \mathcal{N}_g \lra \Omega^\infty \MT{O}{2}$$
and similarly $\alpha_\infty$. The source spaces are connected and so land in a single component: it is not hard to see that it is the component of the Euler characteristic $\chi(N_g)$. However all components of an infinite loop space are homotopy equivalent, so we can translate so that the maps $\alpha_E$ always land in the 0 component of $\Omega^\infty \MT{O}{2}$. Together, Theorem B of \cite{Wahl} and Section 7 of \cite{GMTW} assert that the map
$$\alpha_\infty : B\mathcal{N}_{\infty} \lra \Omega_0^\infty \MT{O}{2}$$
is a homology equivalence.

One can immediately upgrade this theorem: Korkmaz \cite{K} proves that for $g \geq 7$ the index 2 normal subgroup of $\mathcal{N}_{g}$ generated by Dehn twists is perfect, and that the first group homology in this range is $\bZ/2$, so this subgroup is the maximal perfect subgroup. There is then a map
$$\alpha_\infty^+ : B\mathcal{N}^+_\infty \lra \Omega_0^\infty \MT{O}{2}$$
from the plus-construction with respect to this maximal perfect subgroup (obtained by extending $\alpha_\infty$ to the plus-construction, as the map has target an infinite loop space). This is still a homology isomorphism, and $\pi_1(\Omega_0^\infty \MT{O}{2}) \cong H_1(\Omega_0^\infty \MT{O}{2}) \cong \bZ/2$ so it is also a $\pi_1$-isomorphism, so a homotopy equivalence.

The rational cohomology of $\Omega_0^\infty \MT{O}{2}$ is not hard to determine: it is a polynomial algebra on the reductions of certain integral classes $\zeta_i \in H^{4i}(\Omega_0^\infty \MT{O}{2} ; \bZ)$. The author and J. Ebert \cite{ORW} have studied the divisibility of these classes in the integral cohomology of the stable non-orientable mapping class group: they are indivisible.

\subsection{The spectra $\MT{X}{d}$ and tools from stable homotopy theory}
There are two less well-known spectra appearing in this paper, all of the form $\MT{O}{d}$. The first is $\MT{O}{2}$ and has been introduced in the last section. The second is $\MT{O}{1}$\footnote{This spectrum is also known in homotopy theory as $\bR \bP^\infty_{-1}$, though we do not use this notation.} which is constructed in the same way using the unoriented Grassmannians $Gr_1(\mathbb{R}^{n+1})$. A full description of these spectra is available in \cite{GMTW}. We collect here some of their properties.

For a spectrum $\bold{E} = \{ E_n \}$ the associated infinite loop space is
$$\Omega^\infty \bold{E} := \hocolim \Omega^n E_n$$
where the homotopy colimit, or mapping telescope, is taken using the adjoints $E_n \lra \Omega E_{n+1}$ of the structure maps. This construction sends cofibre sequences of spectra to fibre sequences of infinite loop spaces.
By \cite{GMTW} there are cofibre sequences of spectra
\begin{eqnarray}
    \MT{O}{2} &\overset{\omega} \lra \Sigma^\infty BO_{2 +} \overset{\partial} \lra &\MT{O}{1} \label{MTO2Seq} \\ 
    \MT{O}{1} &\overset{\bar{\omega}} \lra \Sigma^\infty BO_{1 +} \overset{\bar{\partial}} \lra &\Sigma^\infty S^0 \label{MTO1Seq}
\end{eqnarray}
and so fibre sequences of associated infinite loop spaces. It should be noted that $\pi_0(\MT{O}{1})=0$, but we will still write $\Omega^\infty_0 \MT{O}{1}$ to remind the reader that this space is connected.

The evaluation maps $\Sigma^n \Omega^n E_n \lra E_n$ give maps on reduced homology
$$\tilde{H}_{*}(\Omega^n E_n) \cong \tilde{H}_{*+n}(\Sigma^n \Omega^n E_n) \lra \tilde{H}_{*+n}(E_n)$$
and taking direct limits on both sides gives a map
$$\sigma_* : \tilde{H}_*(\Omega^\infty \bold{E}) \lra H_*^{spec}(\bold{E})$$
to the spectrum homology of $\bold{E}$. This is known as the \textit{homology suspension}. There is an analogous \textit{cohomology suspension}
$$\sigma^* : H^*_{spec}(\bold{E}) \lra \tilde{H}^*(\Omega^\infty \bold{E}).$$
If $\bold{E}$ is the suspension spectrum of a based space $X$, then $H_*^{spec}(\bold{E}) \cong \tilde{H}_*(X)$ and the inclusion
$$i_* : \tilde{H}_*(X) \lra \tilde{H}_*(Q(X)) = \tilde{H}_*(\Omega^\infty \bold{E})$$
gives a right inverse for $\sigma_*$. Similarly for the cohomology suspension.

For a proper smooth fibre bundle $p: E \lra X$ with fibre $F$ (in other words, $F$ is a compact smooth manifold without boundary), there is a stable map $\trf_p : \Sigma^\infty X_+ \lra \Sigma^\infty E_+$ the \textit{Becker--Gottlieb transfer}. There is a simple description in the case that the base $X$ is compact: then there is an embedding of $E$ into $X \times \bR^n$ over $X$, for some $n$. Choose a tubular neighbourhood $U$ of $E$ in $X \times \bR^n$, which is diffeomorphic to the normal bundle of $E$ in $X \times \bR^n$. Restricted to a fibre $F_x$ embedded in $\bR^n$ this bundle is simply the normal bundle, so the direct sum with the vertical tangent bundle $T_v E$ is trivial. Thus we write $-T_vE$ for the bundle diffeomorphic to $U$. A Pontrjagin--Thom collapse map gives the \textit{pretransfer}
$$prt: \Sigma^n (X_+) \lra U^+ \cong \Th(-T_vE).$$
The inclusion map $-T_vE \lra T_vE \oplus -T_vE$ is proper, so on one-point compactification gives the \textit{inclusion}
$$inc: \Th(-T_vE) \lra \Sigma^n(E_+).$$
The Becker--Gottlieb transfer is the composition $inc \circ prt : \Sigma^n(X_+) \lra \Sigma^n(E_+)$, or its extension to suspension spectra. If the base $X$ is not compact but admits an exhaustion by compact subspaces, then one has a similar definition. In this case one only obtains a map on suspension spectra.

We gather below some results on the Becker--Gottlieb transfer.
\begin{lem}[\cite{BG}, \cite{BM} equation (2.3) on page 137]\label{Transfer}
Let $p$ be a fibre bundle as above. The composition of the induced maps in cohomology $\trf_p^* \circ p^*$ is multiplication by $\chi(F)$, the Euler characteristic of the fibre. Furthermore, if $q: \tilde{E} \lra E$ is another smooth proper fibre bundle, then $p \circ q : \tilde{E} \lra X$ is also a smooth proper fibre bundle and $\trf_{p \circ q}$ and $ \trf_{p} \circ \trf_{q}$ are homotopic.

If $E$ is oriented then the effect of the transfer on cohomology is given by
\begin{equation}\label{Pushforward}
\trf_p^*(x) = p_!(x \cdot e(T_v E))
\end{equation}
where $p_!$ is the Gysin or fibre-integration map corresponding to $p$, and $e(T_v E)$ is the Euler class of the vertical tangent bundle.
\end{lem}

\subsection{Universal definition of characteristic classes}\label{UniversalClasses}
If $c \in H^i(BO_2 ; R)$ is a characteristic class of unoriented 2-plane bundles we can define a characteristic class of unoriented surface bundles:
let $E \overset{\pi} \lra B$ be a surface bundle then $\nu_c(E) := \trf_\pi^*(c(T_vE)) \in H^i(B ; R)$ is a characteristic class, by the naturality properties of the transfer. Classes that arise in this manner have a good interpretation in terms of the spectrum $\MT{O}{2}$.

\begin{defn}
Define a cohomology class $\bar{\nu}_c$ in $\Omega_0^\infty \MT{O}{2}$ by the image of $c$ under
$$H^{*}(BO_2 ; R) \overset{\sigma^*} \lra H^{*}(Q_{0}(BO_{2+}) ; R) \overset{\omega^*} \lra  H^{*}(\Omega_0^\infty \MT{O}{2} ; R).$$
\end{defn}

\begin{thm}
If $\alpha_E : B \lra \Omega_0^\infty \MT{O}{2}$ is the Madsen--Tillmann map for a surface bundle $E \overset{\pi} \lra B$ then $\alpha_E^* (\bar{\nu}_c) = \nu_c(E)$. Thus the $\bar{\nu}_c$ are universal characteristic classes, in the sense that they do not depend on the genus of the fibre of $\pi$.
\end{thm}
\begin{proof}
Let $T: E \to BO_2$ classify $T_v E$, and $\eta: \Th(-T_v E) \lra \MT{O}{2}$ be given by the maps classifying each vertical normal bundle. Then the diagram below commutes.
\begin{diagram}
\Sigma^\infty E_+ & \rTo^{\Sigma^\infty T} & \Sigma^\infty BO_{2+} \\
\uTo^{inc} & & \uTo^{\omega} \\
Th(-T_vE) & \rTo^\eta & \MT{O}{2} \\
\uTo^{prt} & \ruDotsto^{\alpha^\sharp} & \\
\Sigma^\infty B_+ &  & \\
\end{diagram}
By definition $\nu_c = \trf^*(c(T_vE))$, and so it is also given by $(\alpha^\sharp)^* \omega^* (c)$. We wish to identify $\nu_c$ via maps of spaces. Note that the adjoint of $\alpha^\sharp$ is  the map $\alpha_E$. The following diagram commutes by the naturality of the cohomology suspension.
\begin{diagram}
H_{spec}^*(\Sigma^\infty B_+) & = & H^*(B) \\
\uTo^{(\alpha^\sharp)^*} &  & \uTo^{\alpha_E^*}\\
H_{spec}^*(\MT{O}{2}) & \rTo^{\sigma^*} & H^*(\Omega^\infty \MT{O}{2}) \\
\uTo^\omega &  & \uTo^{\Omega^\infty \omega^*} \\
H_{spec}^*(\Sigma^\infty BO_{2+}) & \rTo^{\sigma^*} &  H^*(Q(BO_{2+}))
\end{diagram}
Thus $\nu_c = \alpha_E^* (\Omega^\infty \omega)^* \sigma^*(c) = \alpha_E^*(\bar{\nu}_c)$.
\end{proof}

\begin{example}
The powers of the first Pontrjagin class $p_1^i \in H^{4i}(BO_2 ; \bZ)$ give the characteristic classes $\zeta_i$ defined by Wahl \cite{Wahl}.
\end{example}

\begin{example}
The Stiefel--Whitney classes $w_1^i w_2^j \in H^{i + 2j}(BO_2 ; \bF_2)$ define characteristic classes $\mu_{i, j}$. Even the classes $\mu_{i, 0}$ are difficult to analyse. They satisfy relations $\mu_{i, 0}^2 = \mu_{2i, 0}$, as the cohomology suspension $\sigma^*$ commutes with Steenrod squares, but the algebraic dependence between the odd classes $\mu_{2i+1, 0}$ is not known.
\end{example}

\begin{example}
The same result is true in the oriented case, with $BO_2$ replaced by $BSO_2$. Here the powers of the first Chern class $c_1^i \in H^{2i}(BSO_2 ; \bZ)$ define integral characteristic classes $\kappa_i$ which in the oriented mapping class groups give the Miller--Morita--Mumford classes.
\end{example}

\subsection{Hopf algebras}
We will make much use of the structure of Hopf algebras over a field $\bF$. Let us fix notation, which coincides with that of the classic reference on Hopf algebras of Milnor and Moore \cite{MM}. Recall that an \textit{abelian} Hopf algebra is one that is commutative and cocommutative. All Hopf algebras appearing in this paper are homology groups of infinite loop spaces and as such are abelian. We write $\psi : X \lra X \otimes X$ for the coproduct of a coalgebra $X$, and $\psi_Y : Y \lra X \otimes Y$ for the coaction of the coalgebra $X$ on a right comodule $Y$.

For a coalgebra $A$, a left $A$-comodule $B$ and a right $A$-comodule $C$, the \textit{cotensor product} $B \cotensor_A C$ is defined as the kernel of
$$B \otimes C \overset{\psi_B \otimes id_C - id_B \otimes \psi_C} \lra B \otimes A \otimes C.$$
In general $B \cotensor_A C$ is just a vector space, but if $A$ is cocommutative then it is again an $A$-comodule.

The Hopf algebra \textit{kernel} and \textit{cokernel} of a map $f: A \lra B$ of Hopf algebras are
$A \hker f := A \cotensor_B \bF$ and $B \hcoker f := B \otimes_A \bF$ respectively. When $A$ and $B$ are abelian, kernels and cokernels are again (abelian) Hopf algebras. We will write $PA$ and $QA$ for the vector spaces of primitive and indecomposable elements respectively.

We will make use of two corollaries of Borel's structure theorem for finite-type abelian Hopf algebras over a perfect field (\cite{MM}, Theorem 7.11). 
\begin{lem}\label{Borel}
Firstly, a Hopf subalgebra of a polynomial algebra is again polynomial. Secondly, if $f : A \lra B$ is an injective homomorphism of Hopf algebras, then $B$ is a free $A$-module. Dually, if $f$ is surjective then $A$ is a free $B$-comodule.
\end{lem}

\subsection{The Eilenberg--Moore spectral sequence}\label{EilenbergMoore}
In this section all homology is with coefficients in some field $\bF$. For a fibration sequence $F \lra E \overset{\pi} \lra B$ of connected spaces, the Eilenberg--Moore spectral sequence is
$$E^2_{p,q} = \Cotor_{p, q}^{H_*(B)}(H_*(E) , \bF \,) \Rightarrow H_{p+q}(F)$$
where the subscript on $\Cotor$ denotes the degree $q$ part of the $p$th derived functor of $- \cotensor_{H_*(B)} \bF$. The $H_*(B)$-comodule structure on $H_*(E)$ is induced by the map $\pi_*$. The usual requirement for convergence is that $B$ be simply-connected, but we will use the exotic convergence of \cite{Dwyer}: base spaces will have $\bZ/2$ as fundamental group, which always acts nilpotently on any finitely-generated $\bF_2$-homology group.

When the fibration is one of connected components of infinite loop spaces, 
$$\Omega_0^\infty \bold{F} \lra \Omega_0^\infty \bold{E} \overset{\pi} \lra \Omega_0^\infty \bold{B}$$
coming from a cofibre sequence of finite-type spectra, all the $\bF$-homology groups are finite-type abelian Hopf algebras over $\bF$. Suppose that some condition guaranteeing the convergence of the Eilenberg--Moore spectral sequence holds. If $\pi_* : H_*(\Omega_0^\infty \bold{E} ; \bF \,) \lra H_*(\Omega_0^\infty \bold{B} ; \bF \,)$ is surjective then Lemma \ref{Borel} implies that $H_*(\Omega_0^\infty \bold{E} ; \bF \,)$ is a free $H_*(\Omega_0^\infty \bold{B} ; \bF \,)$-comodule. Thus the $E^2$-term of the Eilenberg--Moore spectral sequence is
$$\Cotor_{p,q}^{H_*(\Omega_0^\infty \bold{B} ; \bF \,)}(H_*(\Omega_0^\infty \bold{E} ; \bF \,) , \bF \,) = H_*(\Omega_0^\infty \bold{E} ; \bF \,) \cotensor_{H_*(\Omega_0^\infty \bold{B} ; \bF \,)} \bF = H_*(\Omega_0^\infty \bold{E} ; \bF \,) \hker \pi_*$$
concentrated in the line $p=0$, so the spectral sequence collapses. Thus
$$\bF \lra H_*(\Omega_0^\infty \bold{F} ; \bF \,) \lra H_*(\Omega_0^\infty \bold{E} ; \bF \,) \overset{\pi_*} \lra H_*(\Omega_0^\infty \bold{B} ; \bF \,) \lra \bF$$
is a short exact sequence of Hopf algebras.

\subsection{Araki--Kudo operations and the homology of $Q(Y_+)$}
An infinite loop space $X$ has a rich structure on its $\bF_2$-homology. Firstly the H-space structure gives the homology the structure of a Hopf algebra, with the Pontrjagin product denoted $*$. Secondly there are Araki--Kudo operations \cite{AK}: these are a family of homomorphisms
$$Q^s : H_i(X ; \bF_2) \lra H_{i+s}(X ; \bF_2)$$
that are natural for maps of infinite loop spaces. These have formal properties similar to the Steenrod squares: a Cartan formula and Ad\'{e}m relations. The reference is \cite{CLM}, and we will use many properties of these operations. We recall some elementary notions. The sequence $I = (s_1,...,s_k)$ is called \textit{admissible} if $s_i \leq 2 s_{i+1}$ for all $i$. The \textit{excess} of such a sequence is the integer $e(I) := s_1 - \sum_{i=2}^k s_i$. The \textit{length} of $I$ is $l(I) := k$. The \textit{degree} of $I$ is $d(I) := \sum_i s_i$. We write $Q^I = Q^{s_1} \cdots Q^{s_k}$. We write $\mathcal{R}$ for the algebra of all operations, and call it the \textit{Dyer--Lashof algebra}.

For a space $Y$, the $\bF_2$-homology of $Q(Y_+)$ has a (non functorial) description as follows. Let $\mathcal{B}$ be a homogenous basis for $H_*(Y ; \bF_2)$. Then the $\bF_2$-homology of $Q(Y_+)$ is the polynomial algebra on the set
$$\{ Q^I (x) \, | \, x \in \mathcal{B}, \, I \text{\, admissible}, \, e(I) > |x| \}.$$
If the space $Y$ is connected, $\pi_0(Q(Y_+)) \cong \bZ$ and we write $[n]$ for the image of $n \in \pi_0$ in $H_0(Q(Y_+))$. The $\bF_2$-homology of the 0 component $Q_0(Y_+)$ is then the polynomial algebra on the set
$$\{ Q^I (x)*[-2^{l(I)}] \, | \, x \in \mathcal{B}, \, I \text{\, admissible}, \, e(I) > |x|, |Q^I(x)| > 0 \}.$$

\section{Mod 2 homology of $\Omega_0^\infty \MT{O}{1}$}\label{MTO1Homology}

It is clear that the map 
$$\bar{\partial} : Q_0(BO_{1 +}) \lra Q_0(S^0)$$
is the transfer for the double covering $EO_1 \lra BO_1$, $t : Q(BO_{1 +}) \lra Q(S^0)$. The effect of this map on homology has been calculated by Mann, Miller and Miller \cite{MMM}:

\begin{lem}
Let $e_i \in H_i(BO_1 ; \bF_2)$ be the unique non-trivial class. Then $t_*(e_i) = Q^i([1])$.
\end{lem}
\begin{cor}\label{MapDeltaBar}
Thus the map $\bar{\partial}: Q_0(BO_{1 +}) \lra Q_0(S^0)$ on homology is 
$$\bar{\partial}_* (e_i * [-1]) = Q^i([1]) * [-2]$$
so is surjective, as these elements generate $H_*(Q_0(S^0))$ over the Dyer--Lashof algebra.
\end{cor}

\begin{proof}[Proof of Theorem \ref{ExactSeqMTO1}]
The Eilenberg--Moore spectral sequence for the fibration
$$\Omega_0^\infty \MT{O}{1} \lra Q_0(BO_{1 +}) \overset{\bar{\partial}} \lra Q_0(S^0)$$
converges to $H_*(\Omega_0^\infty \MT{O}{1})$ as the base space has fundamental group $\bZ/2$. By Corollary \ref{MapDeltaBar}, $\bar{\partial}_*$ is surjective, so by the discussion in Section \ref{EilenbergMoore} we obtain the short exact sequence of Hopf algebras.

Finally, $H_*(Q_0(S^0))$ is a free algebra so the sequence splits as algebras.
\end{proof}

\begin{rem}
If $I$ is some sequence, $Q^I$ is an Araki--Kudo operation and one can iteratedly apply Ad\'{e}m relations to write it as a linear combination of $Q^J$ with $J$ admissible. We write $Q^I = \sum_{J} \lambda_{J}^I Q^J$ where the $\lambda$ are the coefficients needed to express $Q^I$ as admissible monomials. We adopt the convention that we sum over all $J$, and that $\lambda_J^I = 0$ if $J$ is not admissible. Ad\'{e}m relations decrease excess, so $\lambda^I_J = 0$ if $e(J) > e(I)$ also.
\end{rem}

For a sequence $I$ and an integer $i$ such that $I$ is admissible, $e(I) > i$ and $(I, i)$ is not admissible, there is an element $v^{I, i}$ in $QH_*(Q_0(BO_{1 +}))$ given by
$$v^{I, i} := Q^I(e_i)*[-2^{1 + l(I)}] + \sum_{J, j} \lambda^{I, i}_{J, j} Q^J(e_j)*[-2^{1+l(J)}]$$
where the sum is over all sequences $(J, j)$. This clearly lies in the kernel of $Q(\bar{\partial}_*)$. We adopt the convention that $v^{I, i}$ is also defined for $(I, i)$ admissible, and is 0. This is the natural extension, as $Q^I Q^i$ is already written in terms of admissibles.

\begin{thm}\label{MTO1Polynomial}
The Hopf algebra $H_*(\Omega_0^\infty \MT{O}{1} ; \bF_2)$ is polynomial on a set
$$\{ V^{I, i} \,  | \, I \, \text{admissible}, \, e(I) > i , \, (I, i) \, \text{not admissible} \}$$
where $\bar{\omega}_*(V^{I, i}) = v^{I, i}$ modulo decomposables. 
\end{thm}

\begin{proof}It is now clear that $H_*(\Omega_0^\infty \MT{O}{1})$ is a polynomial algebra, as for Hopf algebras over $\bF_2$ a Hopf subalgebra of a polynomial algebra is again polynomial. Furthermore as the sequence splits as algebras there is an exact sequence of indecomposable quotients
$$0 \lra QH_*(\Omega_0^\infty \MT{O}{1}) \overset{Q(\bar{\omega}_*)} \lra QH_*(Q_0(BO_{1 +})) \overset{Q(\bar{\partial}_*)} \lra QH_*(Q_0(S^0)) \lra 0.$$
The elements $v^{I, i}$ described earlier are indecomposable and lie in the kernel of $Q(\bar{\partial}_*)$. Thus there are polynomial generators $V^{I, i}$ such that $Q(\bar{\omega}_*)(V^{I, i}) = v^{I, i}$. The theorem follows from the following Lemma.
\end{proof}

\begin{lem}\label{MTO1PolynomialGenerators}
The $v^{I, i}$ form an additive basis of $\Ker(Q(\bar{\partial}_*))$.
\end{lem}
\begin{proof}
Let $v = \sum \mu_{I, i} Q^I(e_i)*[-2^{1+l(I)}]$ be an element in the kernel of $Q(\bar{\partial}_*)$. Then
\begin{align*}
0 = \bar{\partial}_* (v) &= \sum_{I, i} \mu_{I, i} Q^I Q^i ([1])*[-2^{1+l(I)}] \\
&= \sum_{I, i} \mu_{I, i} \sum_{J, j} \lambda_{J, j}^{I, i} Q^J Q^j ([1])*[-2^{1+l(J)}] \\
&= \sum_{J, j} \left ( \sum_{I, i} \mu_{I, i} \lambda_{J, j}^{I, i} \right ) Q^J Q^j([1])*[-2^{1+l(J)}]
\end{align*}
so $\sum_{I, i} \mu_{I, i} \lambda_{J, j}^{I, i} = 0$ for each $(J, j)$ admissible and of positive excess. Now consider
\begin{align*}
v - \sum_{I, i} \mu_{I, i} v^{I, i} &= \sum_{I, i} \mu_{I, i} \left ( \sum_{J, j} \lambda_{J, j}^{I, i} Q^J (e_j) \right )*[-2^{1+l(J)}] \\
&= \sum_{J, j} \left ( \sum_{I, i} \mu_{I, i} \lambda_{J, j}^{I, i} \right ) Q^J(e_j)*[-2^{1+l(J)}] \\
&= 0.
\end{align*}
So the $v^{I, i}$ span the kernel of the map $Q(\bar{\partial}_*)$.

Now let $\sum_{I, i} \mu_{I, i} v^{I, i} = 0$, so
\begin{align*}
\sum \mu_{I, i} Q^I(e_i)*[-2^{1+l(I)}] &= \sum_{I, i} \mu_{I, i} \sum_{J, j} \lambda_{J, j}^{I, i} Q^J(e_j)*[-2^{1+l(J)}] \\
&= \sum_{J, j} \left ( \sum_{I, i} \mu_{I, i} \lambda_{J, j}^{I, i} \right ) Q^J(e_j)*[-2^{1+l(J)}].
\end{align*}
If $(I, i)$ is not admissible then $Q^I(e_i)*[-2^{1+l(I)}]$ does not appear on the right hand side, so $\mu_{I, i} = 0$. Thus all the coefficients are 0 and the $v^{I, i}$ are linearly independent.
\end{proof}

It will also be important later to understand the action of the Araki--Kudo operations on the polynomial generators $v^{I, i}$ and $V^{I, i}$.

\begin{thm}\label{AKAction}
If $(l, I)$ is an admissible sequence then in $QH_*(Q_0(BO_{2+}))$
$$Q^{l}(v^{I, i}) = v^{(l, I), i} + \sum_{j, J, J'} \lambda_{J, j}^{I, i} \lambda_{J'}^{l, J} v^{J', j}$$
and the same formula holds for $V^{I, i}$ in $QH_*(\Omega_0^\infty \MT{O}{1})$. Furthermore, if $\lambda_{J, j}^{I, i} \neq 0$ then $j > i$.
\end{thm}
\begin{proof}
Firstly we analyse the effect of $Q^l$ on $v^{I, i}*[2^{1+l(I)}]$.
\begin{align*}
Q^l(v^{I, i}*[2^{1+l(I)}]) &= Q^l Q^I(e_i) + \sum_{J, j} \lambda_{J, j}^{I, i} Q^l Q^J(e_j) \\
&= Q^{(l, I)}(e_i) + \sum_{J, j} \lambda_{J, j}^{I, i} \sum_{J'} \lambda_{J'}^{l, J} Q^{J'}(e_{j}).
\end{align*}
Note that this element is in the kernel of $\bar{\partial}_*$, as $v^{I, i}$ is and $\bar{\partial}_*$ commutes with Araki--Kudo operations and the Pontrjagin product. Thus
\begin{align*}
0 &= Q^{(l, I)} Q^i([1]) + \sum_{J', J, j} \lambda_{J, j}^{I, i} \lambda_{J'}^{l, J} Q^{J'} Q^j([1]) \\
&= Q^{(l, I, i)}([1]) + \sum_{J'', J', J, j} \lambda_{J, j}^{I, i} \lambda_{J'}^{l, J} \lambda_{J''}^{J', j} Q^{J''}([1])
\end{align*}
and the terms in the summation are all admissible. Thus this is the unique way to express $Q^{(l, I, i)}$ in terms of admissible operations, so
$$\lambda_{J''}^{l, I, i} = \sum_{J', J, j} \lambda_{J, j}^{I, i} \lambda_{J'}^{l, J} \lambda_{J''}^{J', j}.$$

Secondly we analyse the expression for $v^{(l, I), i}*[2^{2+l(I)}]$.
\begin{align*}
v^{(l, I), i}*[2^{2+l(I)}] &= Q^{(l, I)}(e_i) + \sum_{J''', j'} \lambda_{J''', j'}^{l, I, i} Q^{J'''}(e_{j'}) \\
&= Q^{(l, I)}(e_i) + \sum_{J''', j'} \left ( \sum_{J', J, j} \lambda_{J, j}^{I, i} \lambda_{J'}^{l, J} \lambda_{J''', j'}^{J', j} \right ) Q^{J'''}(e_{j'}) \\
&= Q^{(l, I)}(e_i) + \sum_{J', J, j} \lambda_{J, j}^{I, i} \lambda_{J'}^{l, J} \left ( \sum_{J''', j'} \lambda_{J''', j'}^{J', j} Q^{J'''}(e_{j'}) \right ) \\
&= Q^{(l, I)}(e_i) + \sum_{J', J, j} \lambda_{J, j}^{I, i} \lambda_{J'}^{l, J} \left ( v^{J', j}*[2^{1+l(J')}] + Q^{J'}(e_j) \right ) \\
&= Q^{(l, I)}(e_i) + \sum_{J, j} \lambda_{J, j}^{I, i} Q^{(l, J)}(e_j) + \sum_{J', J, j} \lambda_{J, j}^{I, i} \lambda_{J'}^{l, J} v^{J', j}*[2^{1+l(J')}] \\
&= Q^{l}(v^{I, i}*[2^{1+l(I)}]) + \sum_{J', J, j} \lambda_{J, j}^{I, i} \lambda_{J'}^{l, J} v^{J', j}*[2^{1+l(J')}].
\end{align*}
Then $Q^l(v^{I, i}*[2^{1+l(I)}]) = Q^l(v^{I, i})*[2^{2+l(I)}]$ modulo decomposables, so the result follows after translating back to the 0 component, as $l(J') = l(I)+1$.

Finally, as $(I, i)$ is not admissible but $I$ is, to write $Q^{(I, i)}$ as admissibles an Ad\'{e}m relation must be applied involving the last term $Q^i$. Applying an Ad\'{e}m relation to $Q^a Q^b$ gives monomials $Q^A Q^B$ with $B$ strictly larger than $b$, so if $\lambda_{J, j}^{I, i} \neq 0$ then $j > i$. 
\end{proof}

\section{The map $\partial$ on homology}

We wish to study the composition
$$Q(BO_{2 +}) \overset{\partial} \lra \Omega^\infty \MT{O}{1} \overset{\bar{\omega}} \lra Q(BO_{1 +})$$
and claim that $\partial$ is the pretransfer for the circle bundle $S^1 \lra BO_1 \lra BO_2$. It then follows that $\bar{\omega} \circ \partial$ is the transfer for this bundle.

\begin{rem}
In general there are fibre bundles $S^{n-1} \lra \bS(\gamma_n) \overset{\pi} \lra BO_n$ and a map $f: \bS(\gamma_n) \lra BO_{n-1}$ classifying the vertical tangent bundle. Then the composition
$$Q(BO_{n +}) \overset{\trf_\pi} \lra Q(\bS(\gamma_n)_+) \overset{f}\lra Q(BO_{n-1 +})$$
coincides with the composition
$$Q(BO_{n +}) \overset{\partial} \lra \Omega^\infty \MT{O}{n-1} \overset{\omega} \lra Q(BO_{n-1 +})$$
of maps occurring in the analogs of the fibrations \ref{MTO2Seq} and \ref{MTO1Seq}.
\end{rem}

To evaluate this transfer we will use the technique of Brumfiel--Madsen \cite{BM} of reducing it to a transfer for a finite-sheeted cover. Thus we must find a non-degenerate vector field on the fibre that is equivariant for the action of the structure group. We can not do this directly as $S^1$ does not admit any $O_2$-invariant nondegenerate vector fields, so instead we consider the pullback bundle via $d: BO_1 \times BO_1 \lra BO_2$, the map classifying $\gamma_1 \times \gamma_1$,
$$S^1 \lra \bS(\gamma_1 \times \gamma_1) \lra BO_1 \times BO_1.$$
There is a nondegenerate $O_1 \times O_1$-invariant vector field on $S^1$. It has 4 singular points split into two orbits of two each, having opposite indices. Thus the singular locus $\Sigma \subset \bS(\gamma_1 \times \gamma_1)$ is simply $EO_1 \times BO_1 \coprod BO_1 \times EO_1$. We can thus compute the transfer for this bundle in terms of the transfer map $t$ of Section \ref{MTO1Homology}.

\begin{prop}
The composition $BO_1 \times BO_1 \overset{d} \lra BO_2 \overset{T} \lra Q(BO_{1 +})$ is homotopic to
$$BO_1 \times BO_1 \overset{\Delta} \lra BO_1 \times BO_1 \times BO_1 \times BO_1 \overset{(t \wedge id) \times (id \wedge t)} \lra $$
$$Q_2(EO_1 \times BO_{1 +}) \times Q_2(BO_1 \times EO_{1 +}) \overset{ id \times \chi} \lra Q_2(BO_{1 +}) \times Q_{-2}(BO_{1 +}) \overset{\mu} \lra Q_0(BO_{1 +})$$
\end{prop}
\begin{proof}
This is an elementary application of Theorem 2.10 of \cite{BM}. We are aware that there has been a correction to this paper in \cite{MT}, but our indices are very simple and there is no difference.
\end{proof}

We can now effectively identify the map $\partial_*$. The map $d_*$ is surjective and we proved in the last section that $\bar{\omega}_*$ is injective, so the map $T_* \circ d_*$ determines $\partial_*$.

\begin{thm}\label{PartialFormula}
The composition
$$H_*(BO_1 \times BO_1) \overset{d_*} \lra H_*(BO_{2}) \lra H_*(Q_1(BO_{2 +})) \overset{\bar{\omega}_* \partial_*} \lra H_*(Q_0(BO_{1 +}))$$
is given by
$$e_i \otimes e_j \mapsto \sum_{a=0}^i \sum_{b=0}^j \sum_{s=0}^b \sum_{t=0}^a \binom{b-s}{s} \binom{a-t}{t} Q^{i-a+s}(e_{b-s}) * \chi \left ( Q^{j-b+t}(e_{a-t}) \right )$$
which modulo decomposables is
$$\sum_{s=0}^j \binom{j-s}{s} Q^{i+s}(e_{j-s})*[-2] + \sum_{t=0}^i \binom{i-t}{t} Q^{j+t}(e_{i-t})*[-2].$$
\end{thm}

\begin{proof}
Applying the description of the composition from the previous proposition gives
$$e_i \otimes e_j \mapsto \sum_{a=0}^i \sum_{b=0}^j \left ( Q^{i-a}([1]) \wedge e_b \right ) * \chi \left ( e_a \wedge Q^{j-b}([1]) \right ).$$
We then use the following formulae from \cite{CLM} for evaluating the smash product on homology:
\begin{align*}
  (a * b) \wedge c &= \sum (a \wedge c') * (b \wedge c'') \text{\,\,\,when\,\,\,} \Delta_*(c) = \sum c' \otimes c'' \\  
Q^i(x) \wedge y &= \sum Q^{i+t}(x \wedge Sq^t_*(y)).
\end{align*}
So
$$Q^{i-a}([1]) \wedge e_b = \sum_{s=0}^{b} Q^{i-a+s}([1] \wedge Sq^s_*(e_b))$$
which can be expressed as $\sum_{s=0}^{b} \binom{b-s}{s} Q^{i-a+s}(e_{b-s})$, as $Sq^s_*(e_b) = \binom{b-s}{s} e_{b-s}$. This gives the first expression. The expression modulo decomposables follows immediately, using the fact that $\chi(e_i) = e_i * [-2]$ modulo decomposables.
\end{proof}

\begin{cor}
The formula in Theorem \ref{PartialFormula} for the indecomposable part of the image of $e_i \otimes e_j$ lies in the kernel of $Q(\bar{\partial}_*)$ and so can be expressed as a linear combination of the $v^{I, i}$. This expression is simply
$$ \sum_{s=0}^j \binom{j-s}{s} v^{i+s, j-s} + \sum_{t=0}^i \binom{i-t}{t} v^{j+t, i-t}.$$
Thus $e_0 \otimes e_i$ maps to $v^{i, 0}$ modulo decomposables.
\end{cor}
\begin{proof}
Similarly to the proof of Lemma \ref{MTO1PolynomialGenerators}, consider $$Q(\partial_* \bar{\omega}_* d_*) (e_i \otimes e_j) + \sum_{s=0}^j \binom{j-s}{s} v^{i+s, j-s} + \sum_{t=0}^i \binom{i-t}{t} v^{j+t, i-t}$$ which by the definition of the $v^{I, i}$ can be written as
$$\sum_{s=0}^j \binom{j-s}{s} \sum_{a, b} \lambda_{a, b}^{i+s, j-s} Q^a(e_b)*[-2] + \sum_{t=0}^i \binom{i-t}{t} \sum_{a, b} \lambda_{a, b}^{j+t, i-t} Q^a(e_b)*[-2]$$
so the coefficient of $Q^a(e_b)*[-2]$ is $\sum_{s=0}^j \binom{j-s}{s} \lambda_{a, b}^{i+s, j-s} + \sum_{t=0}^i \binom{i-t}{t} \lambda_{a, b}^{j+t, i-t}$. This element lies in the kernel of $\bar{\partial}_*$, and the $Q^a(e_b)$ are sent to $Q^{(a, b)}([1])*[-4]$ where $(a, b)$ is admissible. These are linearly independent in the homology of $Q_0(S^0)$ and so all the coefficients are 0. Thus the expression above is identically 0, and the result follows.
\end{proof}

\begin{prop}\label{PartialSurjective}
The map $Q(\partial_*) : QH_*(Q_0(BO_{2 +})) \lra QH_*(\Omega_0^\infty \MT{O}{1})$ is surjective.
\end{prop}
\begin{proof}
Define an increasing filtration $G^i := \langle \,V^{a, b} \,\, | \,\, b \leq i \, \rangle$ of $G^\infty := \langle \, V^{a, b} \, \rangle \subset QH_*(\Omega_0^\infty \MT{O}{1})$. The previous corollary implies that $G^0$ is in the image of $Q(\partial_*)$. Consider the indecomposable element $V^{i-a, a}$, so $i-a > 2a$, and $a < i/3$. Then
$$Q(\partial_* d_*) (e_a \otimes e_{i-a}) = \sum_{s=0}^{i-a} \binom{i-a-s}{s} V^{a+s, i-a-s} + \sum_{t=0}^a \binom{a-t}{t} V^{i-a+t, a-t}.$$
Either $2s + a > i$, so $s > i - a - s$ and the binomial coefficient in the first sum is 0, or $2s + a \leq i$. This together with $a < i/3$ implies that $a+s \leq 2(i-a-s)$, so $(a+s, i-a-s)$ is admissible and $V^{a+s, i-a-s}$ is 0. Thus the first sum is 0. The second sum is 
$$V^{i-a, a} + \sum_{t=1}^a \binom{a-t}{t} V^{i-a+t, a-t}$$
so is $V^{i-a, a}$ modulo $G^{a-1}$. Thus by induction along the filtration $G^*$, $G^{\infty} = \langle V^{a, b} \rangle$ is in the image of $Q(\partial_*)$. The map $\partial$ is an infinite loop map, so commutes with Araki--Kudo operations: thus $\mathcal{R} \cdot G^{\infty}$ is also in the image.

Now we introduce a new filtration $F^i := \langle \, V^{I, j} \,\, | \,\, j \geq i \, \rangle$ of $QH_*(\Omega_0^\infty \MT{O}{1})$. This filtration is decreasing, with $F^0$ the whole vector space, but it has the property that in any fixed degree $n$ the filtration $F^*_n$ is has finite length. In particular
$$F^n_n := \langle \, V^{I, j} \,\, | \,\, j \geq n, \, d(I) + j=n \, \rangle = \{ 0 \}$$
as $d(I) > 0$. We will proceed by induction up the filtration in each degree. By Theorem \ref{AKAction}, for $(l, I)$ admissible,
$$Q^{l}(V^{I, i}) = V^{(l, I), i} + \sum_{j, J, J'} \lambda_{J, j}^{I, i} \lambda_{J'}^{l, J} V^{J', j}$$
and for each term in the sum $j > i$. In particular the sum lies in $F^{i+1}$, so $Q^l(V^{I, i}) = V^{(l, I), i}$ modulo $F^{i+1}$. Thus by iterating,
$$V^{(I, k), i} = Q^I(V^{k, i}) \,\,\, \text{modulo} \,\,\, F^{i+1}$$
and so $F^i \subseteq \mathcal{R} \cdot G^{\infty} + F^{i+1}$. In degree $n$, $F^n_n = \{ 0 \}$ is in the image of $Q(\partial_*)$, so by this inclusion $F^{n-1}_n$ is too, and so on: thus $F^{0}_n$ is in the image. This holds in all degrees $n$ so $F^0$ is in the image, but this is the entire space $QH_*(\Omega_0^\infty \MT{O}{1})$.
\end{proof}

\section{The homology of $\Omega_0^\infty \MT{O}{2}$}\label{MTO2Homology}
\begin{proof}[Proof of Theorem \ref{ExactSeqMTO2}]
We now study the Eilenberg--Moore spectral sequence for the fibration
$$\Omega_0^\infty \MT{O}{2} \lra Q_0(BO_{2 +}) \overset{\partial} \lra \Omega_0^\infty \MT{O}{1}$$
noting that $\pi_1(\Omega_0^\infty \MT{O}{1}) \cong \bZ/2$ and so by the discussion in Section \ref{EilenbergMoore} the spectral sequence converges to $H_*(\Omega_0^\infty \MT{O}{2})$. By Proposition \ref{PartialSurjective}, $\partial_*$ is surjective, so by the discussion in Section \ref{EilenbergMoore} we obtain the short exact sequence of Hopf algebras.

Finally, $H_*(\Omega_0^\infty \MT{O}{1})$ is a polynomial algebra, so free, so the sequence is split as algebras.
\end{proof}

\begin{cor}
It now follows that $H_*(\Omega_0^\infty \MT{O}{1} ; \bF_2)$ is also dual to a polynomial algebra, as $H^*(\Omega_0^\infty \MT{O}{1} ; \bF_2)$ injects into $H^*(Q_0 (BO_{2+}) ; \bF_2)$, which is polynomial.
\end{cor}

\subsection{Remarks on integral (co)homology}\label{IntegralDiscussion}
Write $X \left [ \frac{1}{2} \right ]$ for the $\bZ [ \frac{1}{2} ]$-localisation in spaces or spectra. The spectrum $\MT{O}{1} \left [ \frac{1}{2} \right ]$ is contractible (this follows immediately from Corollary 6.3 of \cite{Wahl}), so
$$\MT{O}{2} \left [ \frac{1}{2} \right ] \overset{\simeq}\lra \Sigma^\infty BO_{2 +} \left [ \frac{1}{2} \right ]$$
is an equivalence. There is a zig-zag
$$BO_2 \lra BSO_3 \longleftarrow BSU_2 = \mathbb{H} \mathbb{P}^\infty$$
of maps that induce homology equivalences with $\mathbb{Z} [ \frac{1}{2} ]$-coefficients, so a zig-zag of homotopy-equivalences after $\mathbb{Z} [ \frac{1}{2} ]$-localisation. Thus
$$\Omega^\infty \MT{O}{2} \left [ \frac{1}{2} \right ] \simeq Q(\mathbb{H} \mathbb{P}^\infty_+) \left [ \frac{1}{2} \right ]$$
and on homology $H_*(\Omega_0^\infty \MT{O}{2} ; \bZ [\frac{1}{2}]) \cong H_*(Q_0(\mathbb{H} \mathbb{P}^\infty_+) ; \bZ [\frac{1}{2}])$ as Hopf algebras. The homology of $Q_0(\mathbb{H \, P}^\infty_+)$ can be completely calculated integrally by piecing together its $\bF_p$ homologies, as its Bockstein spectral sequence is a functor of that of $\mathbb{H} \mathbb{P}^\infty$. This determines the integral homology of $\Omega^\infty_0 \MT{O}{2}$ except for its 2-torsion. By Theorem \ref{ExactSeqMTO2} the map $\Omega_0^\infty \MT{O}{2} \lra Q_0(BO_{2 +})$ gives an injection on the $E_1$ pages of their Bockstein spectral sequences, but it is not clear what happens on subsequent pages.

In the range $* \leq 6$ only the primes 2 and 3 contribute. The Bockstein spectral sequence for $Q_0(\mathbb{H \, P}^\infty_+)$ at the prime 3 collapses at the $E_2$ page, so the only odd primary contribution is a $\bZ/3$ in degree 3. In this range one can deduce the Bockstein spectral sequence for $\Omega_0^\infty \MT{O}{2}$ at the prime 2, as the map above still gives an injection on the $E_2$ page, and it collapses at the $E_3$ page. The only unusual contribution is a $\bZ/4$ in degree 3.

\section{A polynomial family in the mod 2 cohomology of the stable non-orientable mapping class group}\label{PolynomialFamily}

Let $F \lra E \overset{\pi} \lra B$ be a bundle of surfaces. We can define mod 2 characteristic classes as follows. There is an associated first real cohomology bundle $V$, with fibre $H^1(F_b ; \bR)$ over the point $b \in B$. Define
$$\chi_i(E) := w_i(V) \in H^i(B ; \bF_2)$$
the $i^{th}$ Stiefel--Whitney class of this vector bundle. The main result is that in the stable non-orientable mapping class group these account for one indecomposable generator in each dimension, and are stably independent.
\begin{thm}\label{ChiClasses}
In the non-orientable mapping class groups, the map
$$\bF_2 [ \chi_1, \chi_2, \chi_3, ... ] \lra H^*(B\mathcal{N}_g ; \bF_2)$$
is an injection in the stable range $* \leq (g-3)/4$.
\end{thm}

The theorem of Korkmaz given in Section \ref{MTODef} implies that $H^1 (B \mathcal{N}_g ; \bF_2) \cong \bF_2$ for $g \geq 7$, so $\chi_1$ is the generator. This gives another interpretation of $\chi_1$: it is the obstruction to reducing the structure group of a bundle of non-orientable surfaces to the index 2 subgroup of the mapping class group generated by Dehn twists.

The proof of the theorem is somewhat indirect. There is a class $1 \in KO^0(E)$ representing the trivial 1-dimensional bundle, and applying the Becker--Gottlieb transfer in real K-theory we obtain a virtual bundle $\trf_\pi^*(1) \in KO^0(B)$. An application of the Atiyah--Singer index theorem for families due to Becker and Schultz implies the following theorem, which gives a homotopy-theoretic characterisation of the virtual bundle $V$.
\begin{thm}
Suppose the surface bundle $E \overset{\pi} \lra B$ is smooth and the fibres are compact, connected and non-orientable. Then the K-theory class $\trf_\pi^*(1)$ coincides with $1 - V$ in $KO^0(B)$. In particular $\chi_i(E) = w_i(-\trf_\pi^*(1))$ in $H^i(B ; \bF_2)$.
\end{thm}
\begin{proof}
By \cite[Theorem 6.1]{BS}, $\trf^*_\pi(1) = \sum_i (-1)^i [H^i(F_b ; \bR)]$. In our situation this sum is $1 - V$ as $[H^0(F_b ; \bR)]$ is the constant rank 1 vector bundle, $[H^1(F_b ; \bR)]$ is the vector bundle $V$, and the higher terms are 0 as the fibres $F_b$ are non-orientable surfaces.
\end{proof}
 
The virtual bundle $-\trf_\pi^*(1)$ is classified by
$$B \overset{\trf} \lra Q(E_+) \overset{collapse} \lra Q(S^0) \overset{Q(i)} \lra \bZ \times BO \overset{\chi} \lra \bZ \times BO$$
where $Q(i)$ is the extension to the free infinite loop space of the inclusion $i$ of $S^0$ to the 0 and 1 components of $\bZ \times BO$, and $\chi$ is the inversion map on $\bZ \times BO$. The discussion in Section \ref{UniversalClasses} implies that the above composition is homotopic to
$$B \overset{\alpha_E} \lra \Omega^\infty \MT{O}{2} \overset{\omega} \lra Q(BO_{2+}) \overset{collapse} \lra Q(S^0) \overset{Q(i)} \lra \bZ \times BO  \overset{\chi} \lra \bZ \times BO$$
as the collapse map $E \lra *$ can be taken to factor through the map $E \lra BO_2$ classifying the vertical tangent bundle. Define $\bar{\chi}_i \in H^i(\Omega_0^\infty \MT{O}{2} ; \bF_2)$ to be the pullback of the $i^{th}$ Stiefel--Whitney class $w_i$ by the composition. We call $\bar{\chi_i}$ the \textit{universal $\chi_i$ class}, as by the previous theorem $\alpha_E^*(\bar{\chi}_i) = \chi_i(E) \in H^i(B ; \bF_2)$ for any smooth, non-orientable surface bundle with compact fibres.

\begin{lem}
The composition $Q(i) \circ \text{collapse} \circ \omega : \Omega_0^\infty \MT{O}{2} \lra BO$ is injective on $\bF_2$-cohomology. The same is true after applying the inversion map $\chi$, as it is a homotopy equivalence.
\end{lem}
\begin{proof}
It is enough to show that the map is surjective on $\bF_2$-homology, and to do this it is enough to show that the composition 
$$\Omega^\infty \MT{O}{2} \overset{\omega} \lra Q(BO_{2+}) \overset{\text{collapse}} \lra Q(S^0) \overset{Q(i)} \lra \bZ \times BO$$
without restriction to 0 components is surjective on $\bF_2$-homology. The action of the Dyer--Lashof algebra on the $\bF_2$-homology of $\bZ \times BO$ has been computed by Priddy \cite{Priddy}, and it is generated over the Dyer--Lashof algebra (and Pontrjagin product) by the class $x_0 \in H_0(\bZ \times BO ; \bF_2)$ representing the component $\{ 1 \} \times BO$. As the composition is a map of infinite loop spaces, it is enough to show that the class $x_0$ is in the image.

The composition $\pi_0(B) \overset{\pi_0(\alpha_E)} \lra \pi_0(\Omega^\infty \MT{O}{2}) \cong \bZ \lra \pi_0(Q(S^0)) \cong \bZ$ for a surface bundle $E \overset{\pi} \to B$ picks out the Euler characteristic of the fibre. In particular, any bundle with fibre $\bR \bP^2$ lands in the 1 component, and so the top composition is an isomorphism on $\pi_0$. In particular, the class $x_0$ is in the image.
\end{proof}
The lemma implies that
$$\bF_2 [\bar{\chi}_1, \bar{\chi}_2, \bar{\chi}_3,...] \lra H^*(\Omega_0^\infty \MT{O}{2} ; \bF_2)$$
is injective, as the Stiefel--Whitney classes are algebraically independent in $H^*(BO ; \bF_2)$. Theorem \ref{ChiClasses} follows by the homology stability of Section \ref{HomStab}.

\bibliographystyle{hplain}
\bibliography{StableHomologyNMCG}

\end{document}